\documentclass[graybox]{svmult}

\usepackage{mathptmx}       
\usepackage{helvet}         
\usepackage{courier}        
\usepackage{type1cm}        
\usepackage{makeidx}         
\usepackage{graphicx}        
\usepackage{multicol}        
\usepackage[bottom]{footmisc}

\makeindex             

\usepackage{amsmath,amssymb,mathtools}

\newcommand{\Z}{\mathbb{Z}}
\newcommand{\R}{\mathbb{R}\ts}

\newcommand{\Q}{\mathbb{Q}\ts}
\newcommand{\N}{\mathbb{N}}

\newcommand{\QQ}{\mathbb{Q}\ts}
\newcommand{\NN}{\mathbb{N}}

\newcommand{\ii}{\ts\mathrm{i}\ts}

\newcommand{\nn}{\nonumber}

\newcommand{\oo}{{\scriptstyle \mathcal{O}}}
\newcommand{\OO}{\mathcal{O}}

\begin{document}

\title*{Well-rounded sublattices and coincidence site lattices}
\author{P. Zeiner}
\institute{P. Zeiner \at Faculty of Mathematics, Bielefeld University, 33501 Bielefeld, Germany, \\ \email{pzeiner@math.uni-bielefeld.de}}
%
%
\maketitle

\abstract{A lattice is called well-rounded, if its lattice vectors of minimal length span the ambient space. We show that there are interesting connections
between the existence of well-rounded sublattices and coincidence site lattices
(CSLs). Furthermore,
we count the number of well-rounded sublattices for several planar lattices
and give their asymptotic behaviour.}

\section{Introduction}

A lattice in $\R^d$ is called well-rounded, if its (non-zero) lattice vectors
of minimal length span $\R^d$. This means that there exist at least $2d$
lattice vectors of minimal positive length, and $\R^d$ has a basis consisting
of lattice vectors of minimal length. However, such a basis need not be
a primitive lattice basis in dimensions $d\geq 4$.

Well-rounded lattices are important for several reasons. Many important lattices
occurring in mathematics and physics are well-rounded. For instance, the
hexagonal lattice and the square lattice in $\R^2$ and the cubic lattices in
$\R^3$ are well-rounded, as are the hypercubic lattices and the $A_4$-lattice
in $\R^4$, which play an important role in quasicrystallography. Examples
in higher dimensions are the Leech lattice, the Barnes-Wall
lattices, and the Coxeter-Todd lattice; see \cite{ConSlo} for background.

Let us briefly mention two problems of mathematical crystallography where
well-rounded lattices occur.
They are connected to the question of densest lattice sphere
packings, as all extreme lattices (those lattices corresponding to densest
lattice sphere packings) are perfect (i.e. the lattice vectors of minimal
length determine the Gram matrix uniquely) and are thus well-rounded.
They also play an important role in reduction theory, as they are exactly
those lattices for which all the successive minima are equal~\cite{vdW}. 

Here, we want to deal with two specific questions: Has a given lattice
well-rounded sublattices, and if so, what are the well-rounded sublattices
and how many are there. The first question is answered in Sec.~\ref{sec:csl}
for planar lattices and a partial answer is given for $d>2$. The second
question is much more difficult in general. Thus we restrict the discussion
to $2$ dimensions, and present some results in Sec.~\ref{sec:plan}.

\section{Well-rounded lattices and CSLs}\label{sec:csl}


Here, we want to deal with the question whether a lattice has a well-rounded
sublattice. It turns out that this question is related to the theory of
coincidence site lattices (CSLs), so let us review the notion of CSL
first. Let $\Lambda$ be a lattice in $\R^d$ and let $R\in O(d)$ be an isometry.
Then $\Lambda(R)=\Lambda\cap R\Lambda$ is called a
\emph{coincidence site lattice (CSL)} if $\Lambda(R)$ is a sublattice of full
rank in $\Lambda$; the corresponding $R$ is called coincidence isometry.
The corresponding index of $\Lambda(R)$ in $\Lambda$
is called coincidence index $\Sigma_\Lambda(R)$, or $\Sigma(R)$ for short. The
set of all coincidence isometries forms a group, which we call $OC(\Lambda)$,
see~\cite{baake97} for details.

Let us look at the planar case first. Here, any two linearly independent
lattice vectors of minimal (non-zero) length form a basis of $\Lambda$.
Let $\gamma$ be the angle between them. Now a well-rounded lattice
is necessarily a rhombic (centred rectangular) lattice such that
$\frac{\pi}{3}<\gamma<\frac{2\pi}{3}, \gamma\ne \frac{\pi}{2}$ 
or a square (corresponding to
$\gamma=\frac{\pi}{2}$) or a hexagonal lattice (corresponding to
$\gamma=\frac{\pi}{3}$ or $\gamma=\frac{2\pi}{3}$). Thus, its symmetry
group is at least $D_2=2mm$, or in other words, there is at least one
reflection symmetry present. As $\Lambda$ and all of its sublattices
have the same group of coincidence isometries~\cite{baake97}, we can infer
that a lattice possesses a well-rounded sublattice only if it has a coincidence
reflection. As the converse holds as well, we have (compare~\cite{BSZ2})
\begin{theorem}\label{lem:wr-csl}
A planar lattice $\Lambda\in \R^2$ has a well-rounded sublattice if and only if
it has a coincidence reflection. 
\end{theorem}

An alternative criterion tells us that a planar lattice has a well-rounded
sublattice if and only if it has a rhombic or rectangular
sublattice~\cite{kuehn}. The existence of well-rounded sublattices can
also be characterised by the entries of the Gram matrices
of $\Lambda$, see~\cite{kuehn}
and~\cite{BSZ2} for various criteria.

One is tempted to generalise these criteria to $d$ dimensions, by using
orthogonal lattices, the $d$-dimensional analogue of rectangular lattices
and orthorhombic lattices in $3$ dimensions. However, this does not work
since a lattice may be well-rounded without having an orthogonal sublattice.
As an example, consider a rhombohedral lattice in $\R^3$, which in general
does not have an orthorhombic sublattice. Nevertheless, an orthogonal
lattice has well-rounded sublattices, and one even has
\begin{theorem}\label{theo:wr-csl-d}
Let $G$ be the symmetry group of an orthogonal lattice, i.e. a lattice that
is spanned by an orthogonal basis. Then $\Lambda$ has a well-rounded
sublattice if \mbox{$G\subseteq OC(\Lambda)$}.
\end{theorem}

This theorem can be proved by induction. The idea is to show that 
\mbox{$G\subseteq OC(\Lambda)$} implies the existence of an orthogonal
sublattice, which in turn implies the existence of well-rounded sublattices.

However, note that the intuitive idea of choosing a ``body-centred orthogonal''
lattice fails in dimensions $d>4$. For if we construct a lattice
as the linear span of the $2^d$ vectors $\sum_{i=1}^d s^{(j)}_i b_i$,
where the $b_i$ form an orthogonal basis of $\R^d$ and
$s^{(j)}_i\in\{1,-1\}$, then these vectors
do not have minimal lengths as at least one of the vectors $2 b_i$ is shorter.
Nevertheless, a modification of this idea works where we choose a suitable
subset of the vectors $\sum_{i=1}^d s^{(j)}_i b_i$.
In particular, if the basis vectors $b_i$ all have approximately the same length
and $d$ is even, we can construct a well-rounded sublattice as the
linear span of
$\sum_{i=1}^d s^{(j)}_i b_i$, where $j$ runs over all possible solutions of
$\sum_{i=1}^d s^{(j)} \equiv 0 \pmod d$.

An immediate consequence of Theorem~\ref{theo:wr-csl-d} is that every
rational lattice has well-rounded sublattices, as $OC(\Lambda)$ contains
all reflections generated by a lattice vector~\cite{zou}.

\section{Well-rounded sublattices of planar lattices}\label{sec:plan}


We now turn to our second question, i.e., we want to find all well-rounded
sublattices of a given lattice. We concentrate on some planar lattices here.
To begin with, we want to find all well-rounded sublattices of the square
lattice. W.l.o.g we may identify it with $\Z^2\simeq\Z[\ii]$.
The idea now is the following. From the previous section, we
know that a planar lattice is well-rounded if and only if it is a rhombic lattice
with $\frac{\pi}{3}<\gamma<\frac{2\pi}{3}$, a square or a hexagonal lattice.
Now a sublattice of a square lattice cannot be hexagonal, so that we can
exclude the latter case, i.e. we only have to find all rhombic and square
well-rounded sublattices. The latter are just the similar sublattices
of the square lattice, which are well known~\cite{baake+grimm,BSZ2}.
The Dirichlet series generating function of their counting function
reads
\begin{align}\label{eq:Dedekind-square}
\Phi^{}_{\square} (s) \, 
  & = \sum_{n\in\N} \frac{s^{}_{\square} (n)}{n^s}
   = \zeta(2s) \Phi^{\mathsf{pr}}_{\square}(s) =
   \, \zeta^{}_{\Q (\ii)} (s)
   \, = \, L (s,\chi^{}_{-4}) \,  \zeta (s)
\end{align}
where $s^{}_{\square} (n)$ is the number of similar sublattices of
the square lattice with index $n$. Here, $\Phi^{\mathsf{pr}}_{\square}(s)$ is
the generating function of the primitive similar sublattices, $\zeta (s)$
is the Riemann zeta function and $\zeta^{}_{\Q (\ii)} (s)$ is the 
Dedekind zeta function
of the complex number field $\Q (\ii)$.

Hence it remains to find all rhombic well-rounded sublattices. Now each
rhombic sublattice has a rectangular sublattice of index $2$, and it is
well-rounded if and only if $\frac{a}{\sqrt{3}}\leq b \leq a \sqrt{3}$
holds, where $a$ and $b$ are the lengths of the orthogonal basis vectors
of the corresponding rectangular sublattice. Thus we only need to
find all rectangular sublattices satisfying the condition above. In fact,
as all square lattices are similar, it is sufficient to find all
rectangular sublattices whose symmetry axes are parallel to those of the
square lattice, and we finally get~\cite{BSZ2}
\begin{align}
 \Phi_{\mathsf{wr},even}(s)\, & = \, \frac{2}{2^s} \Phi^{\mathsf{pr}}_{\square}(s)
   \sum_{p\in\N}\sum_{p < q < \sqrt{3}p}\frac{1}{p^s q^s}     \label{eq:sq-even}\\
 \Phi_{\mathsf{wr}, odd}(s) \, & = \,
   \frac{2}{1+2^{-s}}\Phi^{\mathsf{pr}}_{\square}(s)
   \sum_{k\in\N}\sum_{k < \ell < \sqrt{3}k+\frac{\sqrt{3}-1}{2}}
   \frac{1}{(2k+1)^s (2\ell+1)^s} \label{eq:sq-odd}
\end{align}
where $\Phi_{\mathsf{wr},even}(s)$ and $\Phi_{\mathsf{wr}, odd}(s)$ are the generating
functions counting the rhombic well-rounded sublattices of even and odd indices,
respectively. Putting everything together we arrive at the following
result~\cite{BSZ2}
\begin{theorem}\label{thm:sq1}
  Let\/ $a^{}_{\square} (n)$ be the number of well-rounded sublattices
  of the square lattice with index\/ $n$, and
  $\Phi^{}_{\square,\mathsf{wr}}(s)=\sum_{n=1}^\infty a^{}_{\square}(n)n^{-s}$ 
  the corresponding Dirichlet series generating
  function. It is given by $\Phi^{}_{\square,\mathsf{wr}}(s) =
  \Phi^{}_{\square} (s) + \Phi^{}_{\mathsf{wr},even} (s) + 
  \Phi^{}_{\mathsf{wr},odd} (s)$ with the functions
  from Eqs~\eqref{eq:Dedekind-square},
  \eqref{eq:sq-even} and \eqref{eq:sq-odd}.

If $s>1$, we have the inequality
\[
    D^{}_{\square} (s) - \Phi^{}_{\square} (s) \, < \, 
    \Phi^{}_{\square,\mathsf{wr}}(s) \, < \,
    D^{}_{\square} (s) + \Phi^{}_{\square} (s) \ts .  
\]
with $\Phi^{}_{\square} (s)$ from Eq.~\eqref{eq:Dedekind-square} and the
function
\[
    D^{}_{\square} (s) \, = \,
    \frac{2 + 2^{s}}{1 + 2^{s}} \ts
    \frac{1-\sqrt{3}^{1-s}}{s-1}
    \frac{L(s,\chi^{}_{-4})}{\zeta (2s)} 
    \, \zeta (s) \zeta (2s-1) \ts , 
\]
 
   As a consequence, the summatory function
   $A^{}_{\square} (x) = \sum_{n\leq x} a^{}_{\square} (n)$
   possesses the asymptotic growth behaviour
\[
    A^{}_{\square} (x) \, = \, \frac{\log (3)}{2\ts \pi} \,
      x \log(x)  + \oo(x\log(x))
\]
    as $x\to\infty$.
\end{theorem}
The lower and upper bounds are obtained by approximating the
sums in Eqs.~\eqref{eq:sq-even} and~\eqref{eq:sq-odd} by integrals
via the Euler summation formula, whereas the statement about the asymptotic
behaviour of $A^{}_{\square} (x)$ follows from Delange's theorem, which relates
the asymptotic behaviour of $A^{}_{\square} (x)$ with the analytic properties
of $\Phi^{}_{\square,\mathsf{wr}}(s)$, in particular with its pole at $s=1$.

In fact, we can get additional information about the asymptotic behaviour
of $A^{}_{\square} (x)$ by applying some methods of analytic number theory,
including Dirichlet's hyperbola method and the above mentioned Euler summation
formula (see e.g.~\cite{Apostol}).
\begin{theorem}\label{thm:sq2}
  Let\/ $a^{}_{\square} (n)$ be the number of well-rounded
  sublattices of the square lattice with index\/ $n$. 
  Then, the summatory function $A^{}_{\square} (x) \, =
   \, \sum_{n\leq x} a^{}_{\square} (n)$ possesses the asymptotic
   growth behaviour
\allowdisplaybreaks{
\begin{align}
  A^{}_{\square} (x) 
    & =  \frac{\log(3)}{3} \frac{L(1,\chi^{}_{-4})}{\zeta(2)} x (\log(x) - 1)
            + c_\square x + \OO(x^{3/4}\log(x)) \nn \\
    & = \frac{\log(3)}{2\pi} x \log(x) 
           + \left( c_\square - \frac{\log(3)}{2\pi} \right) x + \OO(x^{3/4}\log(x))
             \nn
\end{align}}
where
\begin{align}
  c_\square& :=  \frac{L(1,\chi^{}_{-4})}{\zeta(2)}
                 \Biggl( \zeta(2) + \frac{\log(3)}{3}
                     \left( \frac{L'(1,\chi^{}_{-4})}{L(1,\chi^{}_{-4})} + \gamma
                              - 2\frac{\zeta'(2)}{\zeta(2)} \right) \nn\\
           & \quad + \frac{\log(3)}{3}\left(2\gamma - \frac{\log(3)}{4} - 
              \frac{\log(2)}{6} \right)   - \sum_{p=1}^\infty \frac{1}{p}
                     \biggl(\frac{\log(3)}{2} - \sum_{p<q< p\sqrt{3}} \frac{1}{q}
                           \biggr)  \nn\\
          &   \quad   - \frac{4}{3}\sum_{k=0}^\infty \frac{1}{2k+1}\biggl(
                          \frac{1}{4}\log(3)
                             - \sum_{k<\ell< k\sqrt{3}+(\sqrt{3}-1)/2} \frac{1}{2\ell+1} 
                                     \biggr) \Biggr) \nn\\
          & \approx  0.6272237  \nn
\end{align}
is the coefficient of $(s-1)^{-1}$ in the Laurent series of
$\sum_n\frac{a^{}_{\square} (n)}{n^s}$ around $s=1$.
Here, $\gamma$ is the Euler-Mascheroni constant.
\end{theorem}


Similar calculations are also possible for the hexagonal lattice. If 
$a^{}_{\triangle} (n)$ is the number of well-rounded
sublattices of the triangular lattice with index $n$,
then the corresponding Dirichlet series generating function
$\Phi^{}_{\triangle,\mathsf{wr}} (s) = \sum_{n=1}^{\infty} a^{}_{\triangle} (n)n^{-s}$
is given by
\[
\Phi^{}_{\triangle,\mathsf{wr}} (s) = \Phi^{}_{\triangle} (s)
  + \Phi^{}_{\triangle,\mathsf{wr}, even} (s) + \Phi^{}_{\triangle,\mathsf{wr},odd} (s),
\]
where
 \begin{equation}\label{eq:Dedekind-tri}
   \Phi^{}_{\triangle} (s) \, = \, \zeta^{}_{\QQ (\rho)} (s)
   \, = \, L (s,\chi^{}_{-3}) \ts \zeta (s) \ts , 
\end{equation}
is the generating function for the similar sublattices of the
hexagonal lattice and
\begin{align}\label{eq:tri-even}
   \Phi^{}_{\triangle,\mathsf{wr},even}(s) \, & = \, \frac{3}{4^s(1+3^{-s})}
   \sum_{p\in\NN}\ \sum_{p < q < 3 p}\frac{1}{p^s q^s}
   \Phi^{\mathsf{pr}}_{\triangle}(s), \\
\label{eq:tri-odd}
  \Phi_{\triangle,\mathsf{wr}, odd}(s) \, & = \, \frac{3}{1+3^{-s}}
  \sum_{k\in\NN}\ \sum_{k < \ell < 3k+1}
  \frac{1}{(2k+1)^s (2\ell+1)^s}
  \Phi^{\mathsf{pr}}_{\triangle}(s)
\end{align}
are the corresponding Dirichlet series for the number of
rhombic well-rounded sublattices with even and odd indices, respectively.
For the asymptotic behaviour we get~\cite{BSZ2}
\begin{theorem}\label{thm:tri2}
  The summatory function $A^{}_{\triangle} (x) \, =
   \, \sum_{n\leq x} a^{}_{\triangle} (n)$ possesses the asymptotic
   growth behaviour
\begin{align}
  A^{}_{\triangle} (x) 
    & =  \frac{9\log(3)}{16} \frac{L(1,\chi^{}_{-3})}{\zeta(2)} x (\log(x) - 1)
              + c_\triangle x + \OO(x^{3/4}\log(x)) \nn \\
    & = \frac{3\sqrt{3}\log(3)}{8\pi} x (\log(x) - 1)
              + c_\triangle x + \OO(x^{3/4}\log(x))  \nn
\end{align}
where $ c_\triangle  \approx 0.4915036  $
is the coefficient of $(s-1)^{-1}$ in the Laurent series of
\ $\sum_n\frac{a^{}_{\triangle} (n)}{n^s}$ around $s=1$.
\hfill$\Box$
\end{theorem}


In both examples, we have infinitely many coincidence reflections, which results
in a large number of well-rounded sublattices and in an asymptotic growth
behaviour of $x\log(x)$. A similar behaviour is to be expected for all
rational lattices, but so far only weaker results have been
obtained~\cite{fuksh-zeta}.

However, in general we have less coincidence
reflections, and we want to conclude with this case.
In fact, if the lattice is not rational, there are either no
or exactly two coincidence reflections~\cite{zou,BSZ2}, and both of them
have the same coincidence index. It is remarkable that in the latter case
the asymptotic behaviour does not depend on the details of the lattice
but only on the coincidence index of its two coincidence reflections. In 
particular we have~\cite{BSZ2} 
\begin{theorem}
  Let\/ $\Lambda$ be a planar lattice that has exactly two coincidence
  reflections. Let $\Sigma$ be their common coincidence index and
  let\/ $a^{}_\Lambda(n)$ denote the
  number of well-rounded sublattices of $\Lambda$ with index\/ $n$. 
  Then, the summatory function $A^{}_{\Lambda} (x) \, =
  \, \sum_{n\leq x} a^{}_{\Lambda} (n)$ possesses the asymptotic
  growth behaviour $
     A^{}_{\Lambda} (x) = 
         \frac{\log 3}{4\Sigma}x + \OO(\sqrt{x}). 
  $
\end{theorem}

\begin{acknowledgement}
The author thanks M. Baake, and R. Scharlau, 
for fruitful discussions. This work was supported by the
German Research Council (DFG) within the CRC 701.
\end{acknowledgement}

\end{document}